\documentclass[11pt]{article} 

\usepackage[utf8]{inputenc} 

\usepackage{geometry} 
\usepackage{graphicx} 

\usepackage{booktabs} 
\usepackage{array} 
\usepackage{paralist} 
\usepackage{verbatim} 
\usepackage{subfig} 
\usepackage{epsfig}
\usepackage{graphicx}
\usepackage{amsmath}
\usepackage{amsfonts}
\usepackage{amssymb}
\usepackage{amsthm}
\usepackage{longtable}
\usepackage{esint}

\usepackage{fancyhdr} 
\usepackage{sectsty}
\allsectionsfont{\sffamily\mdseries\upshape} 

\usepackage[nottoc,notlof,notlot]{tocbibind} 
\usepackage[titles,subfigure]{tocloft} 

\title{High Order Mimetic Symplectic Methods For Hamiltonian Systems}
\author{Anand Srinivasan, Jos\'{e} E. Castillo, \\
Computational Science Research Center, \\
San Diego State University}
\date{arXiv submission: 7-Oct-2023} 

\begin{document}
\maketitle

\section{Abstract}
Hamiltonian systems are known to conserve the Hamiltonian function, which describes the energy evolution over time. Obtaining a numerical spatio-temporal scheme that accurately preserves the discretized Hamiltonian function is often a challenge. In this paper, the use of high order mimetic spatial schemes is investigated for the numerical solution of Hamiltonian equations. The mimetic operators are based on developing high order discrete analogs of the vector calculus quantities divergence and gradient. The resulting high order operators preserve the properties of their continuum ones, and are therefore said to mimic properties of conservation laws and symmetries. Symplectic fourth order schemes are implemented in this paper for the time integration of Hamiltonian systems. A theoretical framework for the energy preserving nature of the resulting schemes is also presented, followed by numerical examples. 

\section{Introduction}

Structure preservation is a desirable numerical property when it comes to solving partial differential equations (PDEs). Schemes that are structure-preserving (or energy-preserving) lead to accurate long term numerical solutions of spatio-temporal systems. Traditional numerical schemes often fail to preserve the numerical energy, and thus lead to numerical instabilities that are non-physical in their behavior. As an example, it is well known that Hamiltonian systems conserve the Hamiltonian function. However, obtaining spatio-temporal discretizations that accurately mimic this structure-preserving property is often a challenge. 

Structure-preserving numerical methods have gained popularity in recent years \cite{furihata2010discrete, wang2018structure, hong2020structure, duan2022convergence}. A structure-preserving scheme preserves positivity, and conserves mass and momentum discretely. These schemes ensure that qualitative features such as motion invariants or structure of the configuration space are reflected in the simulation and therefore provide accurate long term numerical results \cite{sharma2020review}. Most traditional numerical spatio-temporal schemes start with an initial condition and evolve the solution over time as governed by the differential equation. This approach typically fails to discretely mimic the qualitative properties of the continuous system and can eventually result in spurious numerical errors for long term simulations \cite{shadwick2000structure}. In certain cases, such schemes that preserve the invariants have been shown to possess favorable error propagation and are likely to yield more accurate results when compared to a generic method \cite{de1997accuracy}. 

Symplectic numerical integrators are temporal schemes that preserve the symplectic geometric properties of the flow. These structure-preserving temporal schemes are used for the numerical solutions of Hamiltonian systems. A numerical one-step scheme is said to be symplectic if the phase-pair evolution over time is a canonical transformation, resulting in area-preservation \cite{hairer2006geometric}. These methods preserve the Hamiltonian function, conserve momentum and accurately discretize the long-term numerical behavior of the system \cite{sharma2020review}.  

Accurate temporal and spatial discretizations are necessary when it comes to solving a PDE. Spatial operators that discretely mimic a fundamental continuum property are faithful to the physics of the underlying PDE and thereby aid in obtaining more accurate numerical solutions. Consider the Gauss' divergence theorem, where the flux $\vec{F}$ across a surface $\partial S$ is determined by integrating the divergence of the flux over the volume enclosed by the surface. Spatial discretization schemes that discretely mimic this fundamental property yield accurate numerical results. The spatial discretization operators for finite difference schemes are said to be summation by parts operators when they mimic their continuous counterpart of integration by parts \cite{kreiss1974finite}. The spatial discretization at the interiors of the domain can be obtained using centered difference stencils. Extension of the centered difference stencil presents challenges at the boundaries and often result in a reduction in the order of accuracy \cite{strand1994summation}. The equivalent to integration by parts in higher dimensions is the divergence theorem. Obtaining both an accurate functional estimate that discretely mimics the divergence and a discrete quadrature over a volume that produces a discrete quadrature over the surface is often not possible with classical quadrature rules \cite{hicken2013summation}. 

The mimetic schemes considered in this paper are based on the works of Castillo et al \cite{castillo2013mimetic, castillo2003matrix, corbino2020high}. Mimetic methods are based on developing discrete analogs of tensor calculus identities of divergence and gradient, which are used to accurately discretize continuum models for a wide range of physical processes \cite{castillo2001fourth}. These discrete operators preserve the properties of their continuum ones, and thus allow for the discretization of PDEs to mimic critical properties such as conservation laws and symmetries. The resulting quadratures of these high order mimetic methods are diagonal matrices with positive weights, and thus retain positivity of numerical solutions \cite{srinivasan2023mimetic}. The models of boundary value problems solved using the mimetic operators often produce results with more meaningful physical interpretation \cite{castillo2003comparison}. 


The stability and performance of the third order Runge Kutta (RK) scheme with second and fourth order mimetic operators for the advection equation was investigated in \cite{abouali2013stability}, and extended to the shallow water Hamiltonian system of equations in \cite{abouali2014computational}. The energy-preserving nature of the mimetic operators for anisotropic elliptic equations and the advection equation was investigated in \cite{boada2020high, dumett2022energy}. In \cite{castillo2017high}, a theoretical framework for the energy conservation of the wave equation was presented. The position Verlet and Forest Ruth symplectic schemes with mimetic methods for the acoustic wave equation was investigated in \cite{boada2020high2}. The use of mimetic finite difference methods for the Hamiltonian wave equation and its Hamiltonian-preserving properties was investigated in \cite{da2017mimetic}. 

The current work aims to present a theoretical framework for the energy stability of Hamiltonian systems (and the wave equation in particular) using the Corbino-Castillo \cite{corbino2020high} mimetic operators. The use of fourth order explicit symplectic time integration schemes with fourth order mimetic operators is investigated in the context of energy preservation for Hamiltonian systems. The relaxation RK (RRK) temporal scheme was originally applied to the advection equation using mimetic operators in \cite{srinivasan2022mimetic}. The RRK scheme achieves energy conservation at each step of integration. A drawback of this scheme, however, is the computational cost associated with conserving the energy norm. An objective of the current work is to investigate fourth order explicit symplectic integrators that are computationally efficient. 

This paper is organized as follows: we start with an overview of the mimetic schemes and the energy stability for the wave equation. This is followed by a review of explicit temporal schemes that are fourth order and symplectic. Numerical results that validate the energy-preserving discretizations for Hamiltonian systems, along with the computational cost are presented in section 4. 

\section{Mimetic Discretization Methods}
The central idea behind the construction of the mimetic operators is to find high order approximations for the vector calculus quantities of divergence and gradient, and thereby enabling the discretizations to be more faithful to the underlying physics. Consider a smooth region $\Omega \subset \mathbb{R}^3$  of volume $V$ enclosed by the boundary $\partial \Omega$ over the surface $S$. The extended Gauss' divergence theorem states
\begin{equation}\label{EGDT}
	\iiint_{\Omega}  f \ \nabla \cdot \vec{v}\  \ \mathrm{d}V + \iiint_{\Omega} \vec{v} \cdot \left(\nabla f \right) \ \mathrm{d}V  = \oiint_{\partial{\Omega}} f \vec{v} \cdot \vec{n} \ \mathrm{d}S	
\end{equation}
Here, $\vec{n}$ is the outward normal to the boundary, and $f$ and $\vec{v}$ are smooth scalar and vector functions respectively. 
$\nabla \cdot$ is the divergence operator  div, and $\nabla$ the gradient operator grad. In the one dimensional domain $x \in [0,1]$, this becomes integration by parts:
\begin{equation}
	\int_{0}^{1} f \frac{\mathrm{d} v}{\mathrm{d} x}  \ \mathrm{d}x + \int_{0}^{1} \frac{\mathrm{d} f}{\mathrm{d} x} v \ \mathrm{d}x = v(1)f(1) - v(0)f(0),	\label{eq:20}
\end{equation}
where $\dfrac{\mathrm{d} v}{\mathrm{d} x} $ and  $\dfrac{\mathrm{d} f}{\mathrm{d} x}$  are the equivalents of the div and grad operators respectively. The following conservation laws are obtained by setting either $f \equiv 1$ or $v \equiv 1$:
\begin{equation}\label{ConsLaw}
	\int_{0}^{1} \frac{\mathrm{d} v}{\mathrm{d} x}  \ \mathrm{d}x = v(1) - v(0), \quad \quad  \int_{0}^{1} \frac{\mathrm{d} f}{\mathrm{d} x}  \ \mathrm{d}x = f(1) - f(0)
\end{equation} 

Consider the continuous scalar fields $f$ and $g$ (with discrete versions $\mathbf{f}$ and $\mathbf{g}$) and vector fields $\vec{u}, \vec{v}$ (with discrete versions of $\mathbf{u}$ and $\mathbf{v}$). The continuum inner products for these scalar and vectors are defined as 
\begin{equation}
	\int f g \ \mathrm{d} x  \equiv \left<\mathbf{f, g}  \right>_Q  =  \mathbf{f}^T Q \mathbf{g} , \quad \int \vec{u} \ \vec{v} \ \mathrm{d} x \equiv \left<\mathbf{u, v}  \right>_P = \mathbf{u}^T P \mathbf{v}
\end{equation}
Here, $Q$ and $P$ are diagonal matrices with positive weights. That is, $Q$ and $P$ are the quadrature matrices corresponding to the weighted inner products of the scalar and vector functions. The Castillo-Grone and Corbino-Castillo methods \cite{castillo2003matrix, corbino2020high} obtain high order div and grad discretizations that satisfy 
\begin{equation}\label{GDTmim}
	\big  < \hat{\mathbf{D}} \mathbf{v},\hat{\mathbf{f}} \big>_Q + \big <\mathbf{v}, \mathbf{G} \hat{\mathbf{f}} \big>_P = \big <\mathbf{{B}}\mathbf{v}, \hat{\mathbf{f}} \big> 		
\end{equation}
The above equation is a discrete equivalent of the extended Gauss' divergence theorem \cite{corbino2020high}. The discrete versions of the div and grad operator matrices are denoted as $\mathbf{D}$ and $\mathbf{G}$. 

On a staggered one dimensional grid $x \in [0,1]$ comprising of $N$ elements, the spacing length $h = 1/N$.  $\mathbf{D} \in \mathbb{R}^{N \times (N+1)}$ is defined at the cell centers $\mathcal{C} \in \left( i + 
\dfrac{1}{2} \right)h, 0 \leq i \leq (N-1) $ and $\mathbf{G} \in \mathbb{R}^{(N+1) \times (N+2)}$ is defined at the nodes $\mathcal{N} \in ih, 0 \leq i \leq N$. $\mathcal{\hat{C}} \in \left[0, \mathcal{C}, 1  \right] $ is defined by appending the boundary nodes to the cell centers. The functional scalar values $\mathbf{f} \in \mathbb{R}^{N \times 1}$ are evaluated at the cell centers, and is denoted by $\left[ \mathbf{f}_\frac{1}{2}, \mathbf{f}_\frac{3}{2}, \dots,  \mathbf{f}_ {N-\frac{1}{2}} \right]^T$.  $\hat{\mathbf{f}} \in \mathbb{R}^{(N+2) \times 1}$ includes $\mathbf{f}_0$ and $\mathbf{f}_N$ as the first and last elements appended to $\mathbf{f}$ to obtain $ \hat{\mathbf{f}} = \left[\mathbf{f}_0, \mathbf{f}_\frac{1}{2}, \mathbf{f}_\frac{3}{2}, \dots,  \mathbf{f}_ {N-\frac{1}{2}}, \mathbf{f}_N \right]^T $. The vectors $\mathbf{v} \in \mathbb{R}^{(N+1) \times 1}$ are defined at the nodes and denoted by $\left[\mathbf{v}_0, \mathbf{v}_1, \dots, \mathbf{v}_N   \right]^T$. The divergence operator $\mathbf{D}: \mathcal{N} \rightarrow \mathcal{C}$ is a linear mapping that acts on vectors defined at the nodes. The gradient operator $\mathbf{G}: \mathcal{\hat{C}} \rightarrow \mathcal{N}$ acts on scalars defined at the cell centers and boundaries, and linearly maps it to the nodes. $\hat{\mathbf{D}} \in \mathbb{R}^{(N+2) \times (N+1)}: \mathcal{N} \rightarrow \mathcal{\hat{C}}$ is defined by augmenting zeros on the first and last rows of $\mathbf{D}$ to ensure consistency of the matrix dimensions and to account for the fact that the divergence is defined only in the cell centers. In addition, interpolant operators are necessary to map data from nodes to centers and vice-versa. The interpolant $I_D \in \mathbb{R}^{(N+1) \times (N+2)}: \mathcal{N} \rightarrow \mathcal{\hat{C}}$ is a linear mapping for functional values residing at the nodes to the extended cell centers. Similarly, $I_G \in \mathbb{R}^{(N+2) \times (N+1)}: \mathcal{\hat{C}} \rightarrow \mathcal{N}$ maps data from centers to nodes.    
The structure of $\mathbf{B} \in \mathbb{R}^{(N+2) \times (N+2)}$ has the top-left and bottom-right corners of the matrix corresponding to -1 and 1, and zeros otherwise. The diagonal quadrature matrices are of dimensions $Q \in \mathbb{R}^{(N+2) \times (N+2)}$ and $P \in \mathbb{R}^{(N+1) \times (N+1)}$. 

The high order mimetic div and grad operators are derived by satisfying the conservation laws. They possess even order of accuracy $k$ and satisfy the following discrete equivalent of the conservation laws in (\ref{ConsLaw}):
\begin{equation}
	\mathbf{\hat{D}v} - \vec{v}\ ' = \mathcal{O}(h^k), \quad \quad \mathbf{G f} - f \ ' =\mathcal{O}(h^k)
\end{equation} 
The resulting operators do not exactly satisfy (\ref{GDTmim}). Rather, they satisfy (\ref{GDTmim}) up to an error term of order $h$ \cite{srinivasan2023mimetic}. The resulting mimetic boundary operator is denoted as $\mathbf{\hat{B}}$ and is obtained from
\begin{equation}
	Q\mathbf{\hat{D}} + \mathbf{G}^T P = \mathbf{\hat{B}}
\end{equation}

 The resulting boundary operator is comprised predominantly of zeros, with non-zero terms in the rows corresponding to the boundaries. Although this resulting boundary operator $\mathbf{\hat{B}}$ does not satisfy the mimetic constraint in (\ref{GDTmim}) exactly, it has been demonstrated \cite{castillo2003comparison, srinivasan2023mimetic} that $\mathbf{\hat{B}}$ converges to $\mathbf{B}$ as $h$ goes to zero. A consequence of this observation is that when the discrete vectors $\mathbf{v}$ are identical to zero at the boundaries  (i.e., the term on the right hand side of  (\ref{EGDT}) is identical to zero), the duality relation $Q \mathbf{\hat{D}} = - \mathbf{G}^T P$ holds at the interiors of the grid. 

\subsection{Wave Equation And Energy Stability}
Consider the one dimensional two-way wave equation with Dirichlet boundary conditions
\begin{align}
	\frac{\partial^2 u}{\partial t^2} = \nabla \cdot \nabla u + \mathcal{F}, \quad x \in [0, 1], \quad t \geq 0 \label{WaveEq} \\
	u(x,0) = u_0(x), \quad \frac{\partial u}{\partial t}(x,0) = u_1(x), \\
	 u(0,t) =  u(1,t) = 0	 \label{BC}
\end{align}
The discretized mimetic version can be written as 
\begin{align}
	\mathbf{u}_{tt} &= \mathbf{\hat{D}} \mathbf{Gu} + F, \quad t \in [0, T], \\
	\mathbf{u}(0) &= \mathbf{u}_0, \quad  \mathbf{u}_t(0) = \mathbf{u}_1, \\
	\mathbf{u}(0,t) &= \mathbf{u}(1,t) = \mathbf{0}
\end{align}
Well-posedness of the wave equation is guaranteed by implementing the energy method. The energy-stability of the wave equation is obtained by multiplying (\ref{WaveEq}) by $\dfrac{\partial u}{\partial t}$ and integrating over the spatial domain to obtain 
\begin{equation}\label{EnergyStab}
	\int_{0}^{1} \frac{\partial u}{\partial t} \dfrac{\partial^2 u}{\partial t^2} \ \mathrm{d}x = \int_{0}^{1} \frac{\partial u}{\partial t} \ \nabla \cdot \nabla u \ \mathrm{d}x + \int_{0}^{1} \frac{\partial u}{\partial t} \ \mathcal{F} \ \mathrm{d}x 
\end{equation}
The discrete version of the above identity using the mimetic operators becomes
\begin{equation}
	\left<\mathbf{u}_t, \mathbf{u}_{tt}  \right>_Q = \left<\mathbf{u}_t,  \hat{\mathbf{D}} \mathbf{G} \mathbf{u}  \right>_Q + \left< \mathbf{u}_t, F \right>_Q
\end{equation}
The left hand side of the above equation can be written as $\left<\mathbf{u}_t, \mathbf{u}_{tt}  \right>_Q = \dfrac{1}{2} \dfrac{\partial}{\partial t} \left<\mathbf{u}_t, \mathbf{u}_{t}  \right>_Q$. Using integration by parts, the first term on the right hand side of (\ref{EnergyStab}) can be written as follows:
\begin{equation} \label{IBP}
	\int_{0}^{1} u_t \nabla \cdot \nabla u \ \mathrm{d}x = u_t \nabla u \biggr \rvert_{0}^{1} - \frac{1}{2} \frac{\partial}{\partial t} \int_{0}^{1} \left( \nabla u \right)^2 \mathrm{d}x
\end{equation}
As a result of the identity in (\ref{GDTmim}), the discrete mimetic equivalent for the above condition becomes
\begin{equation}
	\left<\mathbf{u}_t, \hat{\mathbf{D}} \mathbf{G} \mathbf{u}  \right>_Q = \left<  \mathbf{\hat{B}}\mathbf{G} \mathbf{u}, \mathbf{u}_t   \right> - \frac{1}{2} \frac{\partial}{\partial t} \left<\mathbf{G}\mathbf{u}, \mathbf{G}\mathbf{u}  \right>_P
\end{equation}
Substituting in the boundary condition from  (\ref{BC}), we obtain
\begin{equation}
	\frac{1}{2} \frac{\partial}{\partial t} \int_{0}^{1} \left[ (u_t)^2 + (\nabla u)^2 \right] \mathrm{d}x = \int_{0}^{1} u_t \ \mathcal{F} \ \mathrm{d}x 
\end{equation}
and the mimetic discrete version becomes
\begin{equation} \label{energyNorm}
	\frac{1}{2} \frac{\partial}{\partial t} \left[ \left<\mathbf{u}_t, \mathbf{u}_t \right>_Q + \left< \mathbf{G}\mathbf{u}, \mathbf{G}\mathbf{u}   \right>_P \right] = \left< \mathbf{u}_t, F \right>_Q
\end{equation}

In the absence of a source term, the conserved Hamiltonian and its discrete equivalent are denoted by
\begin{equation}
\mathcal{H}(u,v) = \dfrac{1}{2} \displaystyle \int_{0}^{1} \left[ (u_t)^2 + (\nabla u)^2 \right] \mathrm{d}x,  \ \    \mathbf{H}(\mathbf{u,v})  =  \dfrac{1}{2} \left[  \left<\mathbf{u}_t, \mathbf{u}_t \right>_Q + \left< \mathbf{G}\mathbf{u}, \mathbf{G}\mathbf{u}   \right>_P \right]
\end{equation}
The Hamiltonian is invariant with respect to time; that is, $\dfrac{\mathrm{d}}{\mathrm{d}t} \mathcal{H}(u,v) = 0$. The discrete Hamiltonian is also preserved when obtained using the mimetic operators (and implementing the duality relationship noted earlier), as shown below:
\begin{align*}
	\frac{\mathrm{d}}{\mathrm{d}t} \mathbf{H}(\mathbf{u,v}) &= \frac{1}{2} \frac{\mathrm{d}}{\mathrm{d}t} \left[  \left<\mathbf{u}_t, \mathbf{u}_t \right>_Q + \left< \mathbf{G}\mathbf{u}, \mathbf{G}\mathbf{u}   \right>_P \right] \\ \nonumber
	&= \frac{1}{2} \left[  \left<\mathbf{u}_{tt}, \mathbf{u}_t \right>_Q + \left< \mathbf{G}\mathbf{u}_t, \mathbf{G}\mathbf{u}   \right>_P \right]   \\ \nonumber
	&=  \frac{1}{2} \left[ \left<\mathbf{\hat{D}Gu}, \mathbf{u}_t \right>_Q  +  \left< \mathbf{G}\mathbf{u}_t, \mathbf{G}\mathbf{u}   \right>_P  \right]   \\ \nonumber
	 &= \frac{1}{2} \left[ -  \left<\mathbf{Gu}, \mathbf{Gu}_t \right>_P  +  \left< \mathbf{G}\mathbf{u}_t, \mathbf{G}\mathbf{u}   \right>_P  \right]  = 0
\end{align*}

In practice, the presence of the non-zero values of $\mathbf{\hat{B}}$ corresponding to the boundaries results in a discrete Hamiltonian function that is preserved to some numeric precision dictated by the grid size. This is illustrated in the numerical results in section 4. Therefore, it can be stated that the high order mimetic operators for Hamiltonian systems results in a time-invariant discretization for the Hamiltonian function. The mimetic operators can therefore be used for the spatial discretization of the wave equation to obtain high order numerical solutions.  

\subsection{Hamiltonian Systems}
Hamiltonian systems \cite{hairer1993solving} of the form 
\begin{equation}
	p_t = - \frac{\partial \mathcal{H}}{\partial q}(p, q), \quad q_t = \frac{\partial \mathcal{H}}{\partial p}(p, q)
\end{equation}
possess the property of preserving the Hamiltonian function $\mathcal{H} = \mathcal{T}(p) + \mathcal{V}(q)$. In addition, the corresponding flow for Hamiltonian systems is symplectic. They are area (and volume)-preserving, and thus conserve the symplectic 2-form wedge product \footnote{$\mathrm{d}p \  \Lambda  \ \mathrm{d}q (\xi_1, \xi_2) = \mathrm{d}p(\xi_1) \ \mathrm{d}q (\xi_2) - \mathrm{d}p(\xi_2) \ \mathrm{d}q (\xi_1)$. } $\mathrm{d}p \  \Lambda  \ \mathrm{d}q$ \cite{yoshida1990construction}. The phase-pair evolution over time $\left( p(0), q(0) \right) \rightarrow (p(\tau), q(\tau))$ is a canonical transformation.  The temporal schemes that retain this volume-preserving behavior are said to be symplectic. The two-way wave equation  (\ref{WaveEq}) can be represented as two first-order in time system of equations	$u_t = v,  v_t = \nabla \cdot \nabla u $. This is a special form of the Hamiltonian system $p_t = q,  q_t = f(p)$ with the conserved Hamiltonian function $\mathcal{H}(p, q) = \dfrac{1}{2} \displaystyle \int \left( q^2 + (\nabla p)^2 \right) \mathrm{d}V $. The numerical solution is obtained using the mimetic Laplacian operator as
\begin{equation} \label{eq:WaveDiscrete}
	\mathbf{u}_t = \mathbf{v}, \quad \mathbf{v}_t = \mathbf{L} \ \mathbf{u}, \quad \text{where} \quad \mathbf{L} =  \mathbf{\hat{D}}  \mathbf{G}
\end{equation}

Following the methodology noted in  \cite{da2017mimetic}, if we denote $\nabla_{u}$ as the gradient operator with respect to the variable $u$, and $\nabla_{v}$ as the gradient with respect to $v$, then
\begin{align}
	Q^{-1} \ \nabla_v \ \mathbf{H}(\mathbf{u,v}) &= Q^{-1} Q \mathbf{v} = \mathbf{v} \\
	Q^{-1} \ \nabla_u \ \mathbf{H}(\mathbf{u,v}) &= Q^{-1} \mathbf{G}^T P \mathbf{Gu} = - Q^{-1} Q \mathbf{\hat{D}Gu} = - \mathbf{\hat{D}Gu}
\end{align}

The semi-discrete equations can thus be written as a Hamiltonian system of ordinary differential equations (ODEs) as
\begin{equation}
\begin{pmatrix}
\mathbf{u}_t \\
\mathbf{v}_t
\end{pmatrix} =  \mathcal{J} Q^{-1} \ \nabla_{uv} \mathbf{H}(\mathbf{u,v}), 
\end{equation}
where $\mathcal{J}$ is the symplectic matrix $\begin{pmatrix}
0 & I \\
-I & 0
\end{pmatrix}  $ and $\nabla_{uv}$ denotes the gradient with respect to the variables $\mathbf{u}$ and $\mathbf{v}$. It is therefore concluded that the spatial discretization obtained using the high order mimetic operators yields a system of ODEs that retains the Hamiltonian character of the given PDE. The Castillo-Grone \& Corbino-Castillo mimetic operators can therefore be considered as a useful scheme for the spatial discretization of Hamiltonian systems. 

\section{Fourth Order Symplectic Integrators}
In this section, an overview of fourth order temporal schemes is presented. The following explicit temporal schemes are investigated: relaxation Runge Kutta scheme \cite{ketcheson2019relaxation}, Forest-Ruth algorithms \cite{forest1990fourth, omelyan2002optimized}  and  composition methods \cite{hairer2006structure, vandekerckhove2014mimetic}. Although not all of these schemes are staggered in time, they possess symplectic (or energy-preserving) properties and are therefore suited for the solution of Hamiltonian systems with the mimetic spatial discretization schemes. 

\subsection{Relaxation Runge Kutta Scheme}
An $s$-stage Runge Kutta method with real coefficients $a_{s1}, a_{s2}, \dots, b_1, \dots, b_s, c_2, \dots, c_s$ of the Butcher tableau \cite{hairer1993solving} with $n$ denoting the time step discretization is given by 
\begin{align}
	\mathbf{u}^{n+1} &= \mathbf{u}^n + \Delta t \sum_{i=1}^s b_i \; \mathbf{k}_i^n, \quad \text{where} \label{eq:StepUpdate} \\
	\mathbf{k}_{s}^n &= \mathbf{u}^n + \Delta t \sum_{j=1}^{i-1} f(t^n + c_j \; \Delta t, a_{ij} \mathbf{k}_{j}^n), \quad i = 1,2,\dots,s       
\end{align}
The $s \times s$-matrix $M$ with elements $m_{ij} = b_i a_{ij} + b_j a_{ji} - b_i b_j$ that satisfies $M = 0$ results in a symplectic RK scheme \cite{sanz1988runge}. In general, the explicit RK schemes are not symplectic (\cite{hairer1993solving}, II.16) \footnote{An example of an explicit RK scheme that is symplectic is presented in \cite{van2008fourth}. }. While implicit schemes such as the Gauss-Radau's are symplectic, the iterative calculations at each time step of integration tend to be computationally expensive. The focus of this paper is restricted to schemes that are explicit in an attempt to seek computationally effective temporal high order schemes. 

The traditional RK schemes fail to conserve  the numerical energy (as defined by the $L_2$-norm) at each step of numerical integration. Taking the norm 
on both sides of eq. (\ref{eq:StepUpdate}), we obtain
\begin{equation}
	||\mathbf{u}^{n+1}||^2 - ||\mathbf{u}^{n}||^2 = \displaystyle 2\Delta t \sum_{j=1}^s b_j \left<\mathbf{k}_j, f_j \right>  - 2\Delta t^2 \sum_{i,j=1}^s b_i a_{ij} \left<f_j, f_i \right> + {\Delta t^2 \sum_{i,j=1}^s b_i b_j \left<f_j, f_i \right>}
\end{equation}
Ketcheson \cite{ketcheson2019relaxation} notes that for conservative systems, the first term on the right hand side of the above equation is zero. However, the other two terms are non-zero and thus lead to a violation of the energy norm. This is remedied by introducing the relaxation parameter $\gamma$ at each time-step, such that the energy norm is satisfied. The energy-preserving relaxation RK step update becomes 
\begin{equation}
	\mathbf{u}^{n+1} = \mathbf{u}^n + \gamma_n \ \Delta t \sum_{i=1}^s b_i \; \mathbf{k}_i^n, \quad  \label{eq:StepUpdateRRK} 
\end{equation}	
so as to ensure $||\mathbf{u}_{\gamma}^{n+1}|| \leq ||\mathbf{u}^n ||$. 

The discrete energy for the wave equation (\ref{eq:WaveDiscrete}) at time $t=0$ is given by $E_0 = \left< \mathbf{v}, \mathbf{v} \right>_Q + \left< \mathbf{G}\mathbf{u}, \mathbf{G}\mathbf{u} \right>_P$, as obtained using the mimetic gradient operator. The step update for the system of equations is then expressed as 
\begin{equation}\label{eq:WaveStepUpdate}
	\mathbf{u}_{\gamma}^{n+1} = \mathbf{u}^n + \gamma^n \ \Delta t \ \mathbf{d}_{u}^n, \quad \quad \\
	\mathbf{v}_{\gamma}^{n+1} = \mathbf{v}^n + \gamma^n \ \Delta t \ \mathbf{d}_{v}^n
\end{equation}
Taking the norm on both sides of  (\ref{eq:WaveStepUpdate}) and substituting $\displaystyle \mathbf{d}^n = \sum_{i=1}^s b_i \; \mathbf{k}_i^n$, we obtain
\begin{equation}\label{eq:quadratic}
	 E^{n+1} = E^n + (\gamma^n)^2 \ \Delta t^2 \ \mathbf{T} + 2 \gamma^n \ \Delta t \ \mathbf{E} + \mathbf{R},
\end{equation}
where $\mathbf{T} = \left<\mathbf{d}_{v}^n, \mathbf{d}_{v}^n  \right>_Q  + \left< \mathbf{G}\mathbf{d}_{u}^n, \mathbf{G}\mathbf{d}_{u}^n  \right>_P$,  $\mathbf{E} = \left<\mathbf{v}^n, \mathbf{d}_{v}^n  \right>_Q  + \left< \mathbf{u}^n, \mathbf{G}\mathbf{d}_{u}^n  \right>_P$ and $\mathbf{R} = \left<\mathbf{v}^n, \mathbf{v}^n  \right>_Q  + \left< \mathbf{u}^n, \mathbf{u}^n  \right>_Q  $. $\gamma^n$ is calculated at each time-step by solving the quadratic in  (\ref{eq:quadratic}). This technique ensures conservation of the energy norm and thus lends itself for application to the solution of Hamiltonian systems. Moreover, this scheme presents the advantage of being explicit. 

When $\Delta t $ is sufficiently small, $\gamma^n$ is close to unity. In \cite{mitsotakis2021conservative}, the authors note that the evaluation of $\gamma^n$ as the roots of the above quadratic equation is often affected by catastropic cancellation due to floating-point error. Alternatively, they suggest the use of a root-solving algorithm such as the bisection method (with a prescribed convergence tolerance) to evaluate the relaxation parameter. In the current work, the use of both the analytical and the root-solver techniques to calculate $\gamma^n$ are investigated as illustrated in the numerical results section. 

\subsection{Forest-Ruth Algorithms}
The Forest-Ruth (FR) algorithm is based on the method of successive canonical transformations to achieve a fourth order scheme that is symplectic. The general form of the algorithm is written as 
\begin{align}
	\mathbf{u}_{i+1} &= \mathbf{u}_{i} - c_i \ \Delta t \  \mathbf{v}_i, \\ \nonumber
	\mathbf{v}_{i+1} &= \mathbf{v}_i + d_i \ \Delta t \  F(\mathbf{u}_{i+1}), \quad i = 1, \dots, 4
\end{align}
The coefficients are given by $c_1 = c_4 = x + 1/2, c_2 = c_3 = -x$, $d_1 = d_3 = 2x + 1, d_2 = -4x - 1, d_4 = 0$, and $x = \dfrac{2^{1/3} + 2^{-1/3} - 1}{6}$. The scheme comprises of three functional evaluations at each time step of integration. It should be noted here that the second intermediate step in the above iteration results in a negative step size, and thus briefly marches backwards in time.  

The position-extended Forest-Ruth (PEFRL) algorithm \cite{omelyan2002optimized} avoids the backwards-step of the FR-algorithm, and results in four functional evaluations at each step of integration. The method is based on deriving Suzuki-like integrators by explicitly reducing the truncation errors to a minimum. The general form of the PEFRL-scheme is as follows:
\begin{align}
	\mathbf{u}_1 &= \mathbf{u}^n + \Delta t \ \xi \mathbf{v}^n,     & \mathbf{v}_1 = \mathbf{v}^n + \Delta t/2 \ (1 - 2 \lambda) \ F(\mathbf{u}_1), \\ \nonumber
	\mathbf{u}_2 &= \mathbf{u}_1 + \Delta t \ \chi \mathbf{v}_1,  	& \mathbf{v}_2 = \mathbf{v}_1 + \Delta t  \ \lambda \ F(\mathbf{u}_2), \\ \nonumber
	\mathbf{u}_3 &= \mathbf{u}_2 + \Delta t \  (1 - 2(\chi + \xi)) \mathbf{v}_2 ,  	& \mathbf{v}_3 = \mathbf{v}_2 + \Delta t  \ \lambda \ F(\mathbf{u}_3), \\ \nonumber
	\mathbf{u}_4 &= \mathbf{u}_3 + \Delta t \  \chi \mathbf{v}_3 , \\ \nonumber
	\mathbf{v}^{n+1} &= \mathbf{v}_3 + \Delta t/2 \ (1 - 2 \lambda) \ F(\mathbf{u}_4), &	\mathbf{u}^{n+1} = \mathbf{u}_4 + \Delta t \  \xi  \ \mathbf{u}^{n+1} 
\end{align}
The coefficients are given by $\xi = +0.1644986515575760E+00$, $ \lambda = - 0.2094333910398989E-01$ and $ \chi = +0.1235692651138917E+01$. 

\subsection{Composition Methods}
Although the FR and PEFRL schemes are symplectic, they do not possess the property of being staggered in time. Temporally staggered schemes are time-reversible in addition to being symplectic \cite{yoshida1990construction}. An intent of the current work is to investigate temporally-staggered schemes for Hamiltonian systems so as to match those of the mimetic schemes. The second order Leapfrog scheme 
\begin{equation}
	\mathbf{u}^{n+1/2} - \mathbf{u}^{n-1/2} = \Delta t \ \mathbf{v}^{n}, \quad \quad \mathbf{v}^{n+1} - \mathbf{v}^{n} = \Delta t \ F(\mathbf{u}^{n+1/2})
\end{equation}
is staggered in space and time and possesses symplectic properties. However, the solutions for $\mathbf{u}$ and $\mathbf{v}$ are at different instances in time (one being calculated at half-integer time steps, while the other at integer steps) and therefore not advantageous when it comes to the computation of numerical energy. The use of composition methods as a staggered-time scheme for the wave equation is presented in \cite{vandekerckhove2014mimetic, verwer2007time}. An $s$-order composition method  is described by the following recursive scheme:
\begin{align}
	\mathbf{u}_k - \mathbf{u}_{k-1} &= -(\beta_k + \alpha_{k-1}) \ \Delta t \ \mathbf{v}_{k-1} \\ \nonumber
	\mathbf{v}_k - \mathbf{v}_{k-1} &= (\beta_k + \alpha_k) \ \Delta t \ \mathbf{u}_k, \quad k = 1, \dots, s
\end{align}
followed by the final step update
\begin{equation}
	\mathbf{u}^{n+1} - \mathbf{u}_s = - \alpha_s \ \Delta t \ \mathbf{v}_s , \quad \mathbf{v}^{n+1} = \mathbf{v}_s
\end{equation}	
The coefficients for the fourth order composition scheme with $s = 5$ are 
$\beta_1 = \alpha_5 = \dfrac{14 - \sqrt{19}}{108}$, $\beta_2 = \alpha_4 = \dfrac{-23 - 20 \sqrt{19}}{270}$, $\beta_3 = \alpha_3 = \dfrac{1}{5}$, $\beta_4 = \alpha_2 = \dfrac{-2 + 10 \sqrt{19}}{135}$, $\beta_5 = \alpha_1 = \dfrac{146 + 5 \sqrt{19}}{540}$. This scheme requires five functional evaluations at each step. 

The numerical investigation of these fourth order temporal schemes with the fourth order mimetic spatial discretization operators is presented in the following section. 

\section{Numerical Examples}
\subsection{Example 1 - Wave equation}
The numerical solution of the wave equation (\ref{WaveEq}) is obtained using the schemes mentioned in the prior section. The following legends are used in the plots below: Runge Kutta - RK4, relaxation Runge Kutta - RRK, Forest-Ruth - FRuth, position-extended Forest Ruth - PEFRL, fourth order composition method - COMP4, relaxation RK scheme with root-finding algorithm - RRK-root. The wave equation was discretized using the mimetic Laplacian operator as shown in (\ref{eq:WaveDiscrete}). The initial condition \\
$u_0(x) = exp \left( -100(x - 0.5)^2 \right)$ and $v_0(x) = 0$ was used, with the spatial domain of $x \in [-30, 30]$. The numerical integration was performed up to $t = 24$ seconds. A CFL-condition of 0.5 was used in the simulations, and the source term was set to zero. 

\begin{figure}
\centering
	  \includegraphics[width=0.75\textwidth]{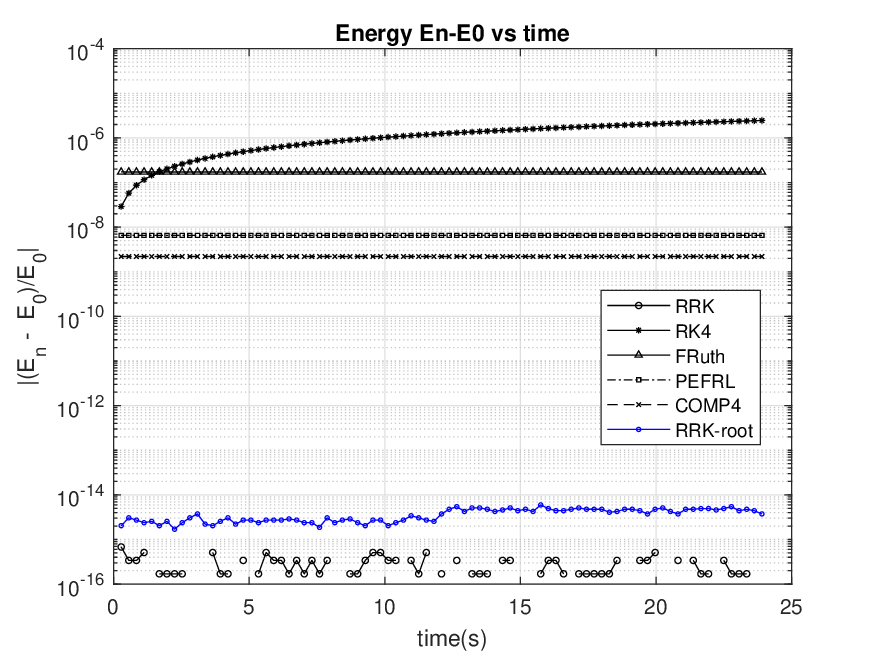}
\caption{Numerical energy vs time, 1D wave equation}
\label{fig:1}       
\end{figure}

Figure (\ref{fig:1}) shows the evolution of the numerical energy for the different schemes. The RK4 scheme shows an asymptotically increasing energy, as noted earlier due to its lack of symplectic-preserving property.  All the other schemes, however, preserve the Hamiltonian. 

\begin{figure}
\centering
	  \includegraphics[width=0.75\textwidth]{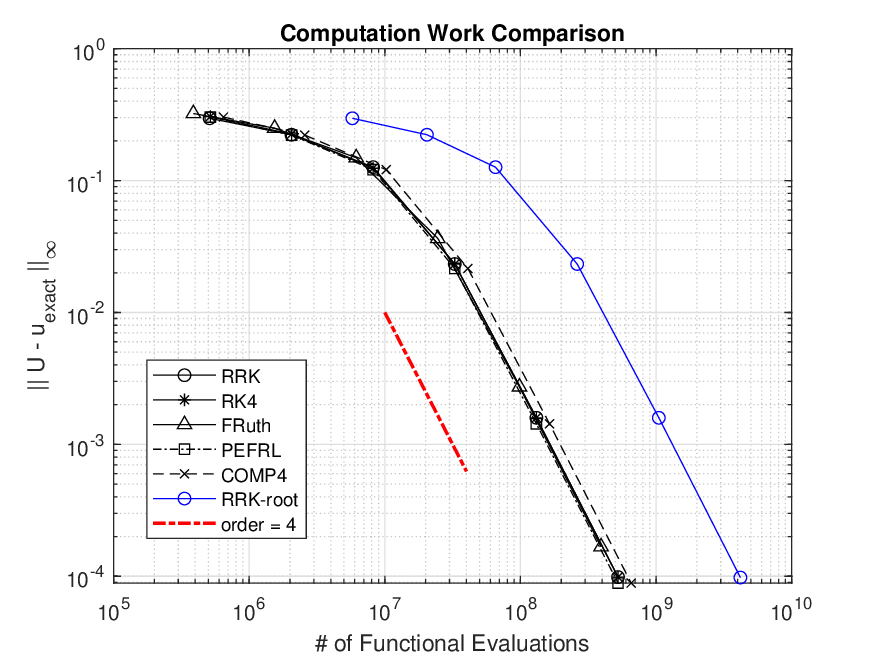}
\caption{Order of convergence, 1D wave equation}
\label{fig:2}       
\end{figure}

Figure (\ref{fig:2}) shows the max-norm of the difference between the exact and numerical solutions as a function of the number of functional evaluations. The schemes exhibit convergence to fourth order as expected. 

\begin{figure}
\centering
	  \includegraphics[width=0.75\textwidth]{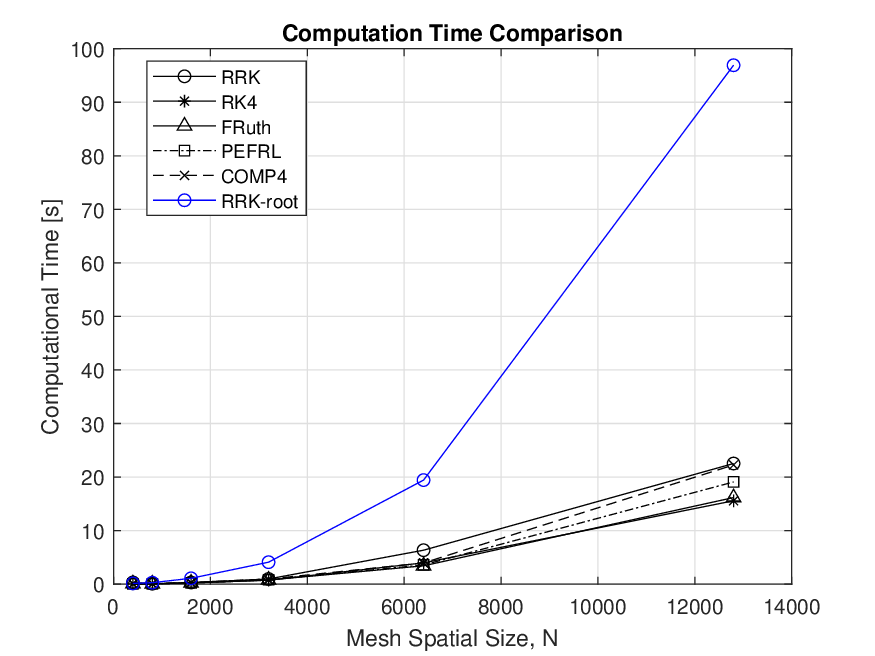}
\caption{Computational time, 1D wave equation}
\label{fig:3}       
\end{figure}

Figure (\ref{fig:3}) shows the computational time comparison for each of the schemes. The RRK-root solver scheme is computationally the most expensive, due to the fact that each iteration requires multiple calls to the bisection root-solver algorithm. The RRK-root solver is about 5-times slower than the other schemes compared in this example. No discernible difference can be noted for the other schemes for the one dimensional problem considered here.  

\subsection{Example 2 - Nonlinear Shallow Water Equations}
\begin{figure}
\centering
	  \includegraphics[width=0.75\textwidth]{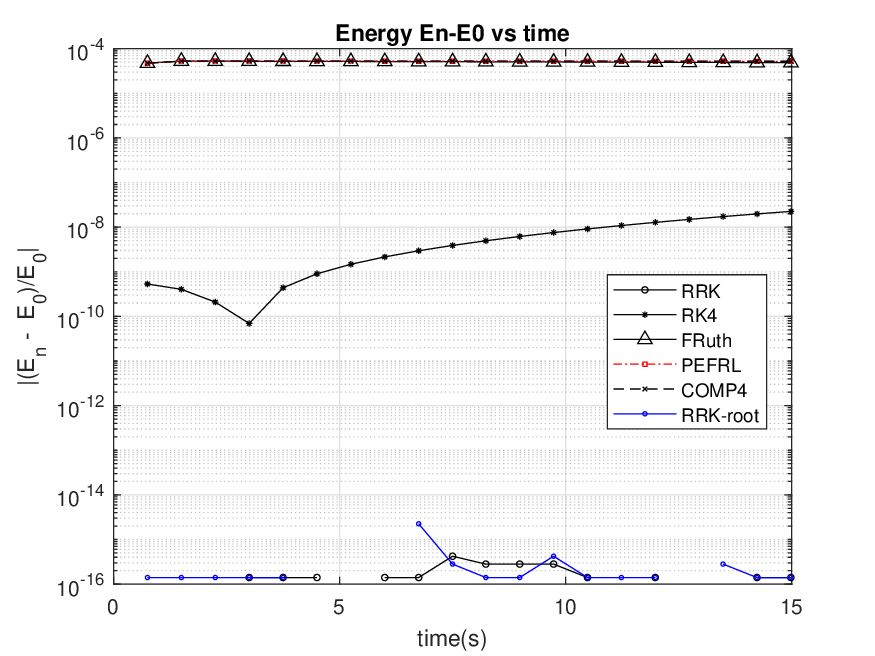}
\caption{1D shallow water equations, energy}
\label{fig:4}       
\end{figure}

\begin{figure}
\centering
	  \includegraphics[width=0.75\textwidth]{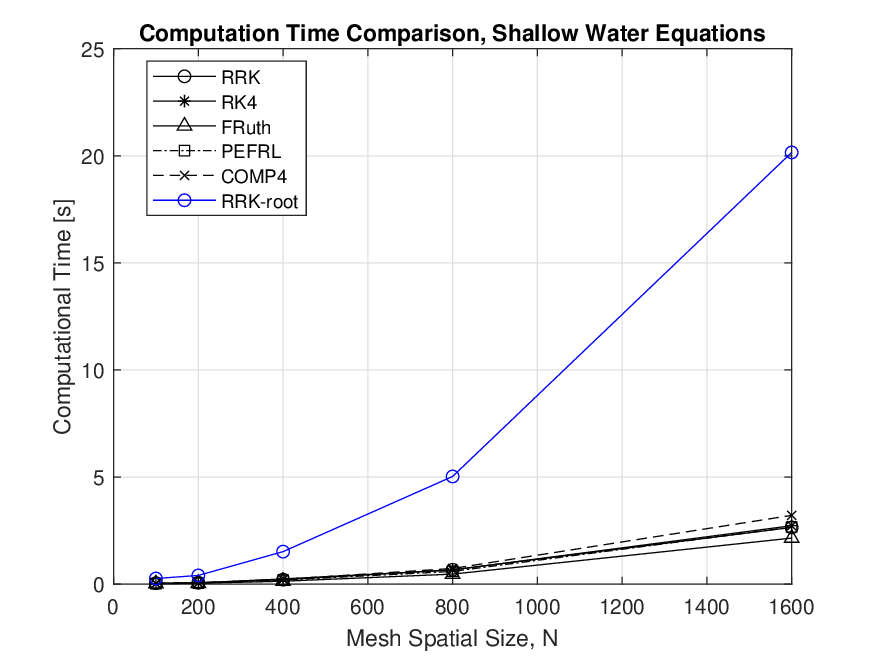}
\caption{1D shallow water equations, computational time}
\label{fig:5}       
\end{figure}

For the second example, the shallow water wave equations 
\begin{equation}
	\eta_t + \nabla \cdot \left[(D + \eta)u  \right] = 0, \quad \quad u_t + g (\ \nabla \cdot u) + (u \cdot \nabla)  u = 0
\end{equation} 
are solved using the schemes noted above. The shallow water equations conserve the Hamiltonian $\mathcal{H}(t, \eta, u) = \dfrac{1}{2} \displaystyle \int_{\Omega} g \eta^2 + (D + \eta)u^2 \ \mathrm{d}x$. The discretized form of the equations are 
\begin{equation}
	\mathbf{e}_t = - \hat{\mathbf{D}}\left[ (D + \mathbf{e}) \mathbf{u} \right] , \quad \mathbf{u}_t = - g \hat{\mathbf{D}} \mathbf{u} - \mathbf{u} \cdot {\mathbf{G}} \mathbf{u}
\end{equation}
The solution was obtained using initial conditions of $\eta(0,x) = 1 + 0.1 exp(-x^2), u(0,x) = 0$ and Dirichlet boundary conditions in the domain $x \in [-30, 30]$, with $D = g = 1$. The resulting energy evolution as a function of time for the various schemes is shown in fig.(\ref{fig:4}) and the computational time comparison is shown in  fig.(\ref{fig:5}). The symplectic schemes conserve the numerical energy for the non-linear system of equations. 

\section{Conclusion}
In this paper, the energy stability of the mimetic discretization methods for Hamiltonian systems is presented. Symplectic time integrators that are explicit temporal schemes have been evaluated with mimetic fourth order operators for Hamiltonian systems. The numerical results show the energy to be discretely preserved for these schemes. 

The RRK schemes computed using an analytical expression for the energy norm show a computational time that is on par with the Forest Ruth and the PEFRL schemes. However, the computational cost increases by about 5-times when solved using a root finding algorithm. Further evaluation is required to assess the computational performance of these schemes in higher dimensions. 

The composition methods are staggered in time and are explicit temporal schemes, thus making it ideally suited for use with mimetic spatial discretization methods for Hamiltonian systems. An advantage of the composition methods is that higher order schemes of 6 exist. It is therefore possible to extend the implementation to higher orders with the corresponding mimetic methods as well.

\bibliography{biblio}
\bibliographystyle{unsrt}

\end{document}